\baselineskip=12pt

\def\adots{\mathinner{\mskip1mu\raise1pt\hbox{.}\mskip2mu\raise4pt\hbox{.}\mskip2mu\raise7pt\vbox{\kern7pt\hbox{.}}\mskip1mu}}

\mathchardef\bfplus="062B
\mathchardef\bfminus="067B
\font\title=cmbx10 scaled\magstep5
\font\chapter=cmbx10 scaled\magstep4
\font\section=cmbx10 scaled\magstep2

\def\~#1{{\accent"7E #1}}
\def\bull{\vrule height .9ex width .8ex depth -.1ex}
\def\sqr#1#2{{\vcenter{\hrule height.#2pt \hbox{\vrule width.#2pt height#1pt \kern#1pt \vrule width.#2pt}\hrule height.#2pt}}}
\def\square{\mathchoice\sqr63\sqr63\sqr{4.2}2\sqr{1.5}2}
\def\Square{\mathop{\square}}
\def\hk#1#2{{\vcenter{\hrule height0.0pt \hbox{\vrule width0.0pt \kern#1pt \vrule width.#2pt height#1pt}\hrule height.#2pt}}}

\vskip 24pt
\centerline {\chapter Real AlphaBeta-Geometries, Walker Geometry,}
\centerline{\chapter and Complex General Relativity}
\vskip 24pt
\noindent {\section Abstract}\hfil\break
By a real $\alpha\beta$-geometry we mean a four-dimensional manifold $M$ equipped with a neutral metric $h$ such that $(M,h)$ admits both an integrable distribution of $\alpha$-planes and an integrable distribution of $\beta$-planes. We obtain a local characterization of the metric when at least one of the distributions is parallel (i.e., is a Walker geometry) and the three-dimensional distribution spanned by the $\alpha$- and $\beta$-distributions is integrable. The case when both distributions are parallel, which has been called two-sided Walker geometry, is obtained as a special case. We also consider real $\alpha\beta$-geometries for which the corresponding spinors are both multiple Weyl principal spinors. All these results have natural analogues in the context of the hyperheavens of complex general relativity.
\vskip 24pt
\noindent $\hbox{Peter R Law}^1$ and $\hbox{Yasuo Matsushita}^2$\hfil\break
${}^1$ 4 Mack Place, Monroe, NY 10950, USA. prldb@member.ams.org\hfil\break
${}^2$Section of Mathematics, School of Engineering, University of Shiga Prefecture, Hikone 522-8533, Japan. matsushita.y@usp.ac.jp\hfil\break
\vskip 24pt
\vfill\eject
\noindent {\section 1. Introduction}
\vskip 12pt
Complex space-times attracted interest in relativity for several reasons including their utility constructing real space-times (Pleba\'nski 1975), Newman's $\cal H$-space (e.g., Newman 1976, Ko et al. 1981), and Penrose's nonlinear graviton (e.g., Penrose 1976, Ward 1980, Penrose and Ward 1980, Penrose 1999). Hyperheavens are complex space-times satisfying the Einstein vacuum equations and with an algebraically degenerate Weyl curvature (Pleba\'nski and Robinson 1976, Finley and Pleba\'nski 1976) and are characterized as solutions of a single partial differential equation called the {\sl hyperheavenly equation}. These various constructions in complex general relativity turned out to be intimately related to each other (e.g, Hansen et al. 1978, Boyer et al. 1980, Newman and Tod 1981).

By a neutral geometry $(M,h)$ we mean a real $2n$-dimensional manifold $M$ equipped with a metric $h$ of neutral signature. One of the intriguing features of neutral geometry is that it manifests analogues with both Riemannian and holomorphic Riemannian geometry, thus mixing themes in definite and indefinite metric geometry. A particularly striking example is the Generalized Goldberg-Sachs Theorem (GGST) in four dimensions (Law 2009, Gover et al. 2011). In the neutral case, a projective spinor field determines a distribution of real totally null planes when real (in complex space-time, a projective spinor field determines a distribution of complex totally null planes) but an almost Hermitian structure when complex (in Riemannian four-manifolds, a projective spinor field determines an almost Hermitian structure). Neutral geometry thus has its own intrinsic appeal but also illuminates both Riemannian geometry and complex general relativity. In particular, the hyperheavenly formalism can be carried over to four-dimensional neutral geometry. The relationship between the geometric structures determined by projective spinor fields and the structure of the Weyl curvature is addressed directly by the GGST but is also prevalent in understanding the geometry of spaces that admit such structures.

We call a neutral geometry admitting a parallel distribution of totally null $n$-planes a Walker geometry, see Walker (1950a). Walker geometry, both in the sense employed here but also in the wider sense of Walker (1950a), has provided a rich source of examples for many geometric topics of interest, see Brozos-V\'azquez et al. (2009). In particular, Walker geometry provides a natural framework in which to formulate the neutral-geometric analogues of results on hyperheavens obtained by Pleba\'nski and co-workers (e.g., Finley and Pleba\'nski 1976 and Boyer et al. 1980) and more generally a context for studying a certain kind of algebraic degeneracy of Weyl curvature. In this paper, all neutral and Walker geometries will be four dimensional.

Walker geometry in four dimensions has a natural spinor description (Law \& Matsushita 2008, to which we refer the reader for a full account of notation, conventions, and details). In particular, a four-dimensional Walker geometry $(M,g)$ has a canonical orientation, with respect to which the parallel distribution is a distribution of self-dual (SD) two-planes, i.e., adopting twistorial notation, a distribution of $\alpha$-planes, which we call an $\alpha$-distribution. Locally, and globally when $(M,g)$ is orientable, the $\alpha$-distribution is equivalent to a projective spinor field $[\pi^{A'}]$: the $\alpha$-plane at $m \in M$ is of the form $\{\,\mu^A\nu^{A'}:\mu^A \in S_m\,\}$, where $S_m$ is the space of unprimed spinors at $m$ (and $S'_m$ that of primed spinors), and $\nu^{A'}$ belongs to the projective class $[\pi^{A'}]$. We denote the $\alpha$-distribution by $Z_{[\pi]}$ and the Walker geometry by $(M,g,Z_{[\pi]})$. The projective spinor field $[\pi^{A'}]$ is a multiple principal spinor of the Weyl curvature. For any (local) projective spinor field $[\pi^{A'}]$, we call a (local) spinor field whose (pointwise) projectivization equals $[\pi^{A'}]$ a {\sl local scaled representative} (LSR) of $[\pi^{A'}]$.

Law (2009) showed how the (real) nonexpanding hyperheavenly equation results readily by imposing the vacuum conditions in a Walker geometry. Law and Matsushita (2011) showed that a general (real) hyperheaven is locally conformally Walker and derived the general (real) hyperheavenly equation from this viewpoint using previously established results on Walker Geometry. More generally, for a neutral geometry with an integrable $\alpha$-distribution $Z_{[\pi]}$, which we call a real $\alpha$-geometry, $[\pi^{A'}]$ is a principal spinor of the Weyl curvature; in Law and Matsushita (2011), we showed that $[\pi^{A'}]$ has multiplicity greater than one iff the neutral geometry is locally conformally Walker, which provides a natural formulation of this part of the GGST.

By a {\sl two-sided} Walker geometry (to be distinguished from a {\sl double} Walker geometry, see Law \& Matsushita (2008a) \S4) is meant a Walker geometry $(M,g,Z_[\pi^{A'}])$ which also admits a parallel distribution of $\beta$ (i.e., anti-self-dual)-planes. This $\beta$-distribution is determined, at least locally, by a projective spinor field $[\lambda_A]$, and so denoted $W_{[\lambda]}$. Chudecki and Przanowski (2008b) exploited the hyperheavenly formalism of Pleba\'nksi and co-workers to provide a local characterization of two-sided Walker geometry. In \S 2, we provide a simple and transparent derivation of this characterization from our perspective in which Walker geometry is treated as fundamental and, in the process, generalize the result to a {\sl sesquiWalker} case by formulating it in the broader context of real $\alpha\beta$-geometry, a four-dimensional neutral geometry $(M,h)$ together with integrable distributions of $\alpha$-and $\beta$-planes. Utilizing earlier work on Walker geometry, we deduce curvature properties for real $\alpha\beta$-geometry and obtain characterizations of certain specializations of curvature (e.g., SD Weyl curvature, left-flatness). In \S 3, we consider algebraically special real $\alpha\beta$-geometry, i.e., when both spinors defined by the null distributions are multiple Weyl principal spinors, and derive conditions for when such are locally conformal to two-sided Walker geometry. We close with an explicit example that illustrates various aspects of our treatment. See Law (2009) for the formalism of null tetrads and spin frames, including spin coefficients and their application to describing null geometry, in the context of neutral geometry.
\vskip 24pt
\noindent {\section 2. Real AlphaBeta-Geometries}
\vskip 12pt
A four-dimensional neutral geometry $(M,h)$ admitting an integrable $\alpha$-distribution will be called a real $\alpha$-geometry and denoted $(M,h,[\pi^{A'}])$, see Law \& Matsushita (2011). The condition for integrability of the $\alpha$-distribution in terms of $[\pi^{A'}]$ is
$$\pi_{A'}\pi^{B'}\nabla_{BB'}\pi^{A'} = 0,\eqno(2.1)$$
where $\pi^{A'}$ is any LSR of $[\pi^{A'}]$. By Law (2009) (6.2.9), any solution $[\pi^{A'}]$  of (2.1) is a Weyl principal spinor (WPS), i.e., a principal spinor (PS)  of the Weyl curvature spinor $\tilde\Psi_{A'B'C'D'}$ (see Law 2006). A Walker geometry $(M,g,[\pi^{A'}])$ is a real $\alpha$-geometry for which $[\pi^{A'}]$ satisfies, in place of (2.1), the stronger condition
$$\pi_{A'}\nabla_{BB'}\pi^{A'} = 0.\eqno(2.2)$$
Equation (2.1) may admit {\sl nontrivially complex} solutions, i.e., solutions $\pi^{A'} \in {\bf C}S'$ which are not just complex scalar multiples of elements of $S'$ (i.e., $\pi^{A'}\bar\pi_{A'} \not= 0$, where $\bar\pi^{A'}$ denotes the complex conjugate of $\pi^{A'}$). Such solutions underly different geometry on $(M,h)$ which we discuss elsewhere. In this paper, we restrict attention to real $\alpha$-geometries and will, for convenience, omit the qualifier `real'.
\vskip 24pt
\noindent {\bf 2.1 Definitions}\hfil\break
An $\alpha\beta$-geometry $(M,h,[\pi^{A'}],[\lambda^A])$ is a neutral geometry admitting both an integrable $\alpha$-distribution $Z_{[\pi]}$ and an integrable $\beta$-distribution $W_{[\lambda]}$. By $\lambda^A$ we denote any LSR of $[\lambda^A]$, and any such LSR satisfies
$$\lambda_A\lambda^B\nabla_{BB'}\lambda^A = 0.\eqno(2.3)$$
Defining the null distributions ${\cal D} := \langle \lambda^A\pi^{A'} \rangle_{\bf R}$ and ${\cal H} := {\cal D}^\perp$, then $(M,h,[\pi^{A'}],[\lambda^A])$ contains the nested null distributions (of types I, II, and III, respectively, in the terminology of Law 2009)
$${\cal D} \leq Z_{[\pi]} \leq {\cal H} \hskip 1.25in {\cal D} \leq W_{[\lambda]} \leq {\cal H}.\eqno(2.4)$$
Note that ${\cal D} = Z_{[\pi]} \cap W_{[\lambda]}$ and ${\cal H} = \langle Z_{[\pi]},W_{[\lambda]} \rangle_{\bf R}$. Equation (2.1) is equivalent to each of:
$$S_b := \pi_{A'}\nabla_b\pi^{A'} = \omega_B\pi_{B'}; \hskip 1in \pi^{B'}\nabla_{BB'}\pi^{A'} =: \eta_B\pi^{A'};\eqno(2.5{\rm a})$$
for some spinors $\omega_B$ and $\eta_B$; while (2.3) is equivalent to each of:
$$\tilde S_b := \lambda_A\nabla_b\lambda^A = \lambda_B\kappa_{B'}; \hskip 1in \lambda^B\nabla_{BB'}\lambda^A =: \zeta_{B'}\lambda^A;\eqno(2.5{\rm b})$$
for some spinors $\kappa_{B'}$, and $\zeta_{B'}$. The significance of the spinors $\omega_A$ and $\eta_{A'}$ was studied in Law (2009), \S 6.2.\bull

A Walker geometry $(M,g,[\pi^{A'}])$ together with an integrable $\beta$-distribution $W_{[\lambda]}$, i.e., an $\alpha\beta$-geometry for which $[\pi^{A'}]$ satisfies (2.2), will be called a {\sl sesquiWalker} $\alpha\beta$-geometry (here the order of $\alpha$ and $\beta$ is meant to indicate which distribution is Walker and which only integrable). A sesquiWalker $\alpha\beta$-geometry for which $[\lambda^A]$ satisfies the analogue of (2.2):
$$\lambda_A\nabla_{BB'}\lambda^A = 0,\eqno(2.6)$$
is of course a two-sided Walker geometry. Note that in the sesquiWalker $\alpha\beta$-case, the choice of the canonical Walker orientation referred to in \S 1 forces the parallel distribution to be an $\alpha$-distribution rather than a $\beta$-distribution; the obvious asymmetry between the $\alpha$- and $\beta$-distributions in a sequiWalker geometry is reflected in the form of the resulting Walker null tetrads and spin frames we employ (defined in Law and Matsushita 2008, \S 2). This asymmetry persists in our {\sl description} of two-sided Walker geometries given below as it is based on the first of the two sets of nested distributions in (2.4).
\vskip 24pt
\noindent {\bf 2.2 Lemma}\hfil\break
For an $\alpha\beta$-geometry $(M,h,[\pi^{A'}],[\lambda^A])$, the distribution $\cal D$  is auto-parallel in the sense of Law (2009) (6.1.7). In a two-sided Walker geometry, both $\cal D$ and $\cal H$ are parallel.

Proof. By (2.5), $\lambda^B\pi^{B'}\nabla_{BB'}\lambda^A\pi^{A'} \propto \lambda^A\pi^{A'}$, i.e., $\cal D$ is auto-parallel. In two-sided Walker geometry, both $Z_{[\pi]}$ and $W_{[\lambda]}$ are parallel, whence ${\cal D} = Z_{[\pi]} \cap W_{[\lambda]}$ and ${\cal H} = {\cal D}^\perp = \langle Z_{[\pi]}, W_{[\lambda]} \rangle_{\bf R}$ are too.\bull
\vskip 24pt
For an $\alpha\beta$-geometry, $Z_{[\pi]}$ and $W_{[\lambda]}$ are each integrable by assumption, and $\cal D$ is integrable being one dimensional. To check the integrability of $\cal H$, it suffices to check whether $[\lambda^B\nu^{B'},\mu^B\pi^{B'}] \in {\cal H}$, for arbitrary $\nu^{B'}$ and $\mu^B$, i.e., whether that expression is orthogonal to $\cal D$. As
$$\eqalign{\lambda_A\pi_{A'}[\lambda^B\nu^{B'},\mu^B\pi^{B'}]^{AA'} &= \lambda_A\pi_{A'}(\lambda^B\nu^{B'}\nabla_{BB'}\mu^A\pi^{A'} - \mu^B\pi^{B'}\nabla_{BB'}\lambda^A\nu^{A'})\cr
&= (\mu^D\lambda_D)(\nu^{D'}\pi_{D'})(\lambda^B\omega_B - \pi^{B'}\kappa_{B'}),\cr}$$
$\cal H$ is integrable iff 
$$\lambda^B\omega_B = \pi^{B'}\kappa_{B'}.\eqno(2.7)$$
One does not, therefore, expect $\cal H$ to be integrable in a general $\alpha\beta$-geometry. In fact, in any (four-dimensional) neutral geometry, the condition for a null distribution $\cal H$ of type III, i.e., ${\cal H}^\perp = {\cal D} := \langle \lambda^A\pi^{A'} \rangle_{\bf R}$, to be integrable is
$$\pi^{A'}\lambda_B\lambda^A\nabla_{AB'}\pi_{A'} = \lambda^A\pi_{B'}\pi^{A'}\nabla_{BA'}\lambda_A,$$ 
see the proof of Law (2009) (6.3.2) where this condition was shown to be equivalent to: $\cal D$ is auto-parallel together with an equation involving spin coefficients. In the context of an $\alpha\beta$-geometry, $\cal D$ is automatically auto-parallel by 2.2, and the previous equation reduces, by (2.5), to (2.7). 
\vskip 24pt
\noindent {\bf 2.3 Lemma}\hfil\break
For a Walker geometry $(M,g,[\pi^{A'}])$ with a $\beta$-distribution $W_{[\lambda]}$, the single condition $\tilde S_b := \lambda_A\nabla_b\lambda^A \propto \lambda_B\pi_{B'}$ is a necessary and sufficient condition for both $W_{[\lambda]}$ and $\cal H$ to be integrable.

For a sesquiWalker $\alpha\beta$-geometry $(M,g,[\pi^{A'}],[\lambda^A])$, the distribution $\cal H$ is integrable iff $\pi^{B'}\kappa_{B'} = 0$, i.e., $\kappa^{A'}$ is an LSR for $[\pi^{A'}]$. For a given LSR $\lambda^A$ of $[\lambda^A]$, there is then an LSR $\pi^{A'}$ of $[\pi^{A'}]$ such that $\tilde S_b = \lambda_B\pi_{B'}$; for a given LSR $\pi^{A'}$ of $[\pi^{A'}]$, there is an LSR $\lambda^A$ of $[\lambda^A]$ such that $\tilde S_b = \lambda_B\pi_{B'}$. Moreover, $\cal H$ is integrable iff auto-parallel. 

We will call a (sesquiWalker) $\alpha\beta$-geometry for which $\cal H$ is integrable, an {\sl integrable (sesquiWalker) $\alpha\beta$-geometry}. Clearly, two-sided Walker geometry is a special case of integrable sesquiWalker $\alpha\beta$-geometry.

Proof. The first assertion is clear from (2.5b) and (2.7). With $Z_{[\pi]}$ parallel, $\cal H$ is auto-parallel iff, for any spinor $\nu^{A'}$, the covariant derivative of $\lambda^A\nu^{A'}$, where $\lambda^A$ is any LSR of $[\lambda^A]$, along $\cal H$ lies in $\cal H$, i.e., is orthogonal to $\cal D$, i.e., $\lambda_A\pi_{A'}\lambda^B\nabla_{BB'}\lambda^A\nu^{A'} = 0$ and $\lambda_A\pi_{A'}\pi^{B'}\nabla_{BB'}\lambda^A\nu^{A'} = 0$. The first equality is equivalent to (2.3) and the second to $\pi^{B'}\kappa_{B'} = 0$. (Of course, any auto-parallel distribution is necessarily integrable, see Law 2009 (6.1.7).)\bull 
\vskip 24pt
The interest in integrability of $\cal H$ is not just passing curiosity.
\vskip 24pt
\noindent {\bf 2.4 Construction}\hfil\break
Beginning with an integrable $\alpha\beta$-geometry $(M,h,[\pi^{A'}],[\lambda^A])$, choose Frobenius coordinates $(p,q,x,y)$ for the first of the nested distributions in (2.4), i.e., so that ${\cal D} = \langle \partial_p \rangle_{\bf R}$, $Z_{[\pi]} = \langle \partial_p,\partial_q \rangle_{\bf R}$, and ${\cal H} = \langle \partial_p,\partial_q,\partial_x \rangle_{\bf R}$. Since $y$ is constant on the integral manifolds of each of the distributions, $dy = \lambda_A\pi_{A'}$ for some LSRs $\lambda_A$ of $[\lambda_A]$ and $\pi_{A'}$ of $[\pi_{A'}]$; as $x$ is constant on the integral surfaces (called $\alpha$-surfaces) of $Z_{[\pi]}$, then $dx = \mu_A\pi_{A'}$, for some spinor $\mu_A$ satisfying $\lambda^A\mu_A \not= 0$. At this stage, one can proceed to construct new coordinates $(u,v,x,y)$ with respect to which the metric takes a coordinate form generalizing the Walker coordinate form, see Law \& Matsushita (2011) (3.19). The construction of these coordinates is essentially the first step in the hyperheavenly formalism. Such coordinates are Frobenius for $Z_{[\pi]}$ but will not, in general, respect the nesting of $\cal D$ within $Z_{[\pi]}$.

We therefore specialize the context to that of an integrable sesquiWalker $\alpha\beta$-geometry $(M,g,[\pi^{A'}],[\lambda^A])$. The assumption that $\cal H$ is integrable is thus characterized by $\pi^{A'}\kappa_{A'} = 0$. With $dy = \lambda_A\pi_{A'}$ and $dx = \mu_A\pi_{A'}$ as in the previous paragraph, one can apply the construction of Law \& Matsushita (2008), 2.3, to obtain Walker coordinates for $(M,g,[\pi^{A'}])$, i.e., one obtains coordinates $(u,v,x,y)$ which are Frobenius coordinates for $Z_{[\pi]}$, yield Walker's canonical coordinate form for the metric
$$\left(g_{\bf ab}\right) = \pmatrix{{\bf 0}_2&{\bf 1}_2\cr {\bf 1}_2&W\cr},\qquad\hbox{where}\qquad W = \pmatrix{a&c\cr c&b\cr}\eqno(2.8)$$
with $a$, $b$, and $c$ some functions of the coordinates, and moreover
$$\partial_u = \mu^A\pi^{A'} \qquad \partial_v = \lambda^A\pi^{A'} \qquad dx = \mu_A\pi_{A'} \qquad dy = \lambda_A\pi_{A'}.\eqno(2.9)$$
As in Law \& Matsushita (2008), 2.4, it proves convenient to specialize the choice of Walker coordinates to {\sl oriented} Walker coordinates. First, by rescaling the LSRs as follows: $\pi^{A'} \mapsto \gamma\pi^{A'}$; $\lambda^A \mapsto \lambda^A/\gamma$; $\mu^A \mapsto \mu^A/\gamma$; one preserves (2.9) but, by appropriate choice of $\gamma$, can ensure $\lambda^A\mu_A = \pm 1$. The LSRs are now fixed up to a common sign. If $\lambda^A\mu_A = 1$, write $\mu_A$ as $\alpha_A$, so that (2.9) becomes
$$\partial_u = \alpha^A\pi^{A'} \qquad \partial_v = \lambda^A\pi^{A'} \qquad dx = \alpha_A\pi_{A'} \qquad dy = \lambda_A\pi_{A'}.\eqno(2.10{\rm a}).$$
In this case, $(v,u,x,y)$ are Frobenius coordinates respecting the nested distributions (note, in particular, that the coordinate tangent vectors of $x$ with respect to the two coordinate systems $(p,q,x,y)$ and $(u,v,x,y)$ differ by an element of $Z_{[\pi]}$; whence, as that with respect to $(p,q,x,y)$ lies in $\cal H$, so does that with respect to $(u,v,x,y)$) and $(u,v,x,y)$ are {\sl oriented} Walker coordinates (i.e., satisfying Law \& Matsushita 2008, (2.8)).

If, however, $\lambda^A\mu_A = -1$, then, as in Law \& Matsushita (2008), 2.4, to achieve oriented Walker coordinates one can resort to Law \& Matsushita (2008), (A1.7), interchanging $u$ with $v$ and $x$ with $y$. After relabelling the coordinates and writing $\mu_A$ as $\beta_A$, one obtains oriented Walker coordinates $(u,v,x,y)$, with (2.9) now taking the form
$$\partial_u = \lambda^A\pi^{A'} \qquad \partial_v = \beta^A\pi^{A'} \qquad dx = \lambda_A\pi_{A'} \qquad dy = \beta_A\pi_{A'},\eqno(2.10{\rm b})$$
and where now $(u,v,y,x)$ are Frobenius coordinates respecting the nested distributions. 

Alternatively, when $\lambda^A\mu_A = -1$, and noting that the Walker Lagrangian in Law \& Matsushita (2008), A1.2, is also invariant under $x \mapsto -x$, $u \mapsto -u$ and $c \mapsto -c$, one can replace $x$ by $-x$ and $u$ by $-u$ (which in effect replaces $\mu_A$ by $-\mu_A$). After relabelling the coordinates one obtains oriented Walker coordinates $(u,v,x,y)$, with (2.9) now taking the same form as (2.10a) but with $\alpha_A = -\mu_A$, and where $(v,u,x,y)$ are again Frobenius coordinates for the nested distributions.

Thus, by an appropriate tactic, one can always find oriented Walker coordinates $(u,v,x,y)$ satisfying (2.10a), where $\lambda^A\alpha_A = 1$, and with $(v,u,x,y)$ Frobenius coordinates respecting the nested distributions. From Law \& Matsushita (2008), (2.11), the Walker spin frames associated to these oriented Walker coordinates are $\{\alpha^A,\lambda^A\}$ and $\{\pi^{A'},\xi^{A'}\}$, i.e., $\beta^A = \lambda^A$ and in terms of the associated Walker null tetrad, one has 
$$\displaylines{{\cal D} = \langle \partial_v \rangle_{\bf R} = \langle \tilde m^a \rangle \qquad Z_{[\pi]} = \langle \partial_u, \partial_v \rangle_{\bf R} = \langle \ell^a, \tilde m^a \rangle_{\bf R} \qquad {\cal H} = \langle \partial_u, \partial_v, \partial_x \rangle_{\bf R} = \langle \ell^a, n^a, \tilde m^a \rangle_{\bf R}\cr
\noalign{\vskip -6pt}
\hfill\llap(2.11{\rm a})\cr
\noalign{\vskip -6pt}
W_{[\lambda]} = \langle \tilde m^a, n^a \rangle_{\bf R} = \langle \partial_v,-{a \over 2}\partial_u - {c \over 2}\partial_v + \partial_x \rangle_{\bf R} = \langle \partial_v,\partial_x - {a \over 2}\partial_u \rangle_{\bf R}.\cr}$$
For the form (2.10b), one has instead $\lambda^A = \alpha^A$ in the Walker spin frames, whence
$$\displaylines{{\cal D} = \langle \partial_u \rangle_{\bf R} = \langle \ell^a \rangle \qquad Z_{[\pi]} = \langle \partial_u, \partial_v \rangle_{\bf R} = \langle \ell^a, \tilde m^a \rangle_{\bf R} \qquad {\cal H} = \langle \partial_u, \partial_v, \partial_y \rangle_{\bf R} = \langle \ell^a, m^a, \tilde m^a \rangle_{\bf R}\cr
\noalign{\vskip -6pt}
\hfill\llap(2.11{\rm b})\cr
\noalign{\vskip -6pt}
W_{[\lambda]} = \langle \ell^a, \tilde m^a \rangle_{\bf R} = \langle \partial_u,{c \over 2}\partial_u + {b \over 2}\partial_v - \partial_x \rangle_{\bf R} = \langle \partial_u,\partial_y - {b \over 2}\partial_v \rangle_{\bf R}.\cr}$$
Note that it is the assumption of integrability of $\cal H$ which allows one to write $dy = \lambda_A\pi_{A'}$ and ultimately achieve the oriented Walker coordinates satisfying (2.10a) (or (2.10b)). Without that assumption, all one can do is construct oriented Walker coordinates without any relationship to $\cal D$ or $\cal H$ (the integrability of $\cal D$ is automatic and does not facilitate matters).\bull 
\vskip 24pt
\noindent {\bf 2.5 Proposition}\hfil\break
Let $(M,g,[\pi^{A'}],[\lambda^A])$ be an integrable sesquiWalker $\alpha\beta$-geometry. Then, for oriented Walker coordinates satisfying the form (2.10a), the metric components in (2.8) satisfy $a_v = 0$, i.e., the component $g_{22}$ is constant along the integral curves of $\cal D$. Note that this coordinate condition is somewhat geometric in nature in that: $\cal D$ is determined by the sesquiWalker $\alpha\beta$-geometry; ${\cal H} = \langle Z_{[\pi]},\partial_x \rangle_{\bf R}$ and $a = g(\partial_x,\partial_x)$.

For the form (2.10b), one obtains instead that $b = g(\partial_y,\partial_y)$ is constant along the integral curves of $\cal D$, i.e., $b_u = 0$, which has an analogous geometric interpretation.

If $(M,g,[\pi^{A'}],[\lambda^A])$ is a 2-sided Walker geometry, then $\cal H$ is integrable and for oriented Walker coordinates satisfying the form (2.10a), $a_v = c_v = 0$, i.e., the metric components $a$ and $c$ are constant along integral curves of $\cal D$. For oriented Walker coordinates satisfying the form (2.10b), $b_u = c_u = 0$, i.e., the metric components $b$ and $c$ are constant along the integral curves of $\cal D$.  Chudecki \& Przanowski (2008b), \S 5, gave one version of this result for 2-sided Walker geometry.

Proof. For any spinor $\kappa^A$ and null tetrad,
$$\nabla_b\kappa^A = n_bD\kappa^A + \ell_bD'\kappa^A - \tilde m_b\delta\kappa^A - m_b\triangle\kappa^A.\eqno(2.12)$$
If $(M,g,[\pi^{A'}],[\lambda^A])$ is an integrable sesquiWalker $\alpha\beta$-geometry, then $[\lambda^A]$ satisfies (2.3). As $\cal H$ is integrable, one can exploit oriented Walker coordinates satisfying (2.10). For form (2.10a), $\lambda_A = \beta_A$ in the Walker spin frame, so by (2.12),
$$0 = \beta_A\beta^B\nabla_b\beta^A = \beta^B\ell_b\beta_AD'\beta^A - \beta^Bm_b\beta_A\triangle\beta^A = \pi_{B'}(\beta_AD'\beta^A) - \xi_{B'}(\beta_A\triangle\beta^A).$$
By Law (2009) (5.8), $\beta_AD'\beta^A = -(a_v/2)$ while $\beta_A\triangle\beta^A = 0$, which proves the relevant assertion. For the form (2.10b), $\lambda^A = \alpha^A$, and the analogous computation yields
$$0 = \alpha_A\alpha^B\nabla_b\alpha^A = -\xi_{B'}\alpha_AD\alpha^A + \pi_{B'}\alpha_A\delta\alpha^A = -{b_u \over 2}\pi_{B'}.$$
In passing, note that, by 2.3, integrability of $\cal H$ can be stated as $\lambda_A\pi^{B'}\nabla_b\lambda^A = 0$. Exploiting the form (2.10a) with $\lambda^A = \beta^A$, one finds from (2.12), $\beta_A\pi^{B'}\nabla_b\beta^A = \beta_A\pi^{B'}n_bD\beta^A - \beta_A\pi^{B'}m_b\triangle\beta^A = 0$, since the Walker spin frames are parallel with respect to $D$ and $\triangle$, i.e., on $\alpha$-surfaces of $Z_{[\pi]}$(Law 2009 (5.8)). This computation confirms that, apart from providing the construction of oriented Walker coordinates satisfying (2.10), the assumption of integrability of $\cal H$ involves no further conditions on the metric components.

Now suppose $(M,g,[\pi^{A'}],[\lambda^A])$ is 2-sided Walker; in particular $[\pi^{A'}]$ satisfies (2.2) and $[\lambda^A]$ (2.6). Hence, in (2.5), $\omega^A$ and $\kappa^{A'}$ are zero and (2.7) is trivially satisfied, i.e., $\cal H$ is integrable. Thus, one can construct oriented Walker coordinates satisfying (2.10). For form (2.10a), with $\lambda^A = \beta^A$, one obtains, using (2.12) and formulae for the action of $D$, $D'$, $\delta$ and $\triangle$ on the elements of the Walker spin frames given in Law (2009) (5.8),
$$0 = \beta_A\nabla_b\beta^A = \beta_A\left[\ell_b\left({a_v \over 2}\alpha^A + {c_v-a_u \over 4}\beta^A\right) - \tilde m_b\left(-{c_v \over 2}\alpha^A + {c_u - b_v \over 4}\beta^A\right)\right] = -{a_v \over 2}\ell_b - {c_v \over 2}\tilde m_b,$$
whence $a_v = c_v = 0$ as claimed. The analogous computation for form (2.10b) with $\lambda^A = \alpha^A$ is
$$0 = \alpha_A\nabla_b\alpha^A = \alpha_A\left[\ell_b\left({a_u - c_v \over 4}\alpha^A + {c_u \over 2}\beta^A\right) - \tilde m_b\left({b_v-c_u \over 4}\alpha^A - {b_u \over 2}\beta^A\right)\right] = {c_u \over 2}\ell_b + {b_u \over 2}\tilde m_b.$$\bull
\vskip 24pt
By Law (2009) (6.2.9), for an $\alpha\beta$-geometry $(M,h,[\pi^{A'}],[\lambda^A])$, $[\pi^{A'}]$ is a WPS of $\tilde\Psi_{A'B'C'D'}$ and $[\lambda^A]$ a WPS of $\Psi_{ABCD}$. By Law \& Matsushita (2008), 2.5, for a sesquiWalker $\alpha\beta$-geometry, $[\pi^{A'}]$ is in fact a multiple WPS and also a PS of $\Phi_{ABA'B'}$, i.e., a {\sl Ricci Principal Spinor (RPS)}. Hence, for a 2-sided Walker geometry, $[\lambda^A]$ is also a multiple WPS and a RPS. These facts are readily seen to be consistent with Proposition 2.5 by consulting the expressions for the components of the curvature spinors with respect to the Walker spin frames given in Law \& Matsushita (2008) \S 2. For example, for form (2.10a), substituting $a_v = c_v = 0$ into Law \& Matsushita (2008) (2.25) yields $\Psi_3 = \Psi_4 = 0$, i.e., $[\beta^A]$ is indeed a multiple WPS. Substituting $a_v = c_v = 0$ into Law \& Matsushita (2008) A1.8, shows that $\nu = 0$, and into A1.7 that $\zeta = 0$ (these quantities being defined at those locations in Law \& Matsushita 2008); hence, in Law \& Matsushita (2008) (2.32--33), one sees that $A_{AB}\beta^A\beta^B = B_{AB}\beta^A\beta^B = 0$, whence $\Phi_{ABA'B'}\beta^A\beta^B = 0$. In an integrable sesquiWalker $\alpha\beta$-geometry, $a_v = 0$ yields only $\Psi_4 = 0$ in Law \& Matsushita (2008) (2.25), as expected. For form (2.10b) in the two-sided Walker case, one obtains instead that $\Psi_0 = \Psi_1 = 0$ from Law \& Matsushita (2008) (2.25) (i.e., $[\alpha^A]$ is a multiple WPS), $\mu = 0$ from Law \& Matsushita (2008) A1.8 and $\Upsilon = 0$ from Law \& Matsushita (2008) A1.7, which entail in Law \& Matsushita (2008) (2.32--33) that $B_{AB}\alpha^A\alpha^B = A_{AB}\alpha^A\alpha^B = 0$.

The covariant derivatives of the Walker null tetrad and spin frames in a Walker geometry were given in Law (2009) (5.5) and (5.8) respectively. For an integrable sesquiWalker $\alpha\beta$-geometry $(M,h,[\pi^{A'}],[\lambda^A])$, using coordinates of either form of (2.10), one readily computes from (2.12) and Law (2009), (5.8), that
$$\lambda^B\nabla_{BB'}\lambda^A \propto \pi_{B'}\lambda^A, \qquad\hbox{whence}\qquad \lambda^B\nabla_{BB'}\lambda^A\pi^{A'} \propto \pi_{B'}\lambda^A\pi^{A'};\eqno(2.13)$$
consequently, from (2.5b), $[\zeta_{B'}] = [\pi_{B'}] = [\kappa_{B'}]$.
For two-sided Walker geometries one has in addition 
$$D'\lambda^A \propto \lambda^A,\qquad \delta\lambda^A \propto \lambda^A,\qquad\nabla_{BB'}\lambda^A\pi^{A'} \propto \lambda_B\pi_{B'}\lambda^A\pi^{A'}.\eqno(2.14)$$
From Law \& Matsushita (2008) (A1.8), in integrable sesquiWalker geometry, the Walker coordinate parametrising integral curves of $\cal D$ ($v$ for form (2.10a), $u$ for form (2.10b)) is an affine parameter for such curves as null geodesics. These geodesics are the intersections of $\alpha$- and $\beta$-surfaces, and the tangent vector ($\partial_v = \lambda^A\pi^{A'}$ for (2.10a), $\partial_u = \lambda^A\pi^{A'}$ for (2.10b)) along such a null geodesic is in fact parallel over the $\alpha$-surface in which the null geodesic lies, and, in the two-sided Walker case, parallel over both the $\alpha$- and $\beta$-surface whose intersection is that null geodesic. The final equation in (2.14) implies that $\cal D$ is parallel, as it must be in two-sided Walker geometry because both $Z_{[\pi]}$ and $W_{[\lambda]}$ are parallel, and of course ${\cal H} = {\cal D}^\perp$ is also parallel.

The null distributions $\langle \partial_u \rangle_{\bf R}$ and $\langle \partial_v \rangle_{\bf R}$ for Walker coordinates were studied in Law (2009) (6.1.33) and (6.1.44) and those results therefore provide a local description of the distribution $\cal D$ here. That discussion confirms that $\cal D$ is auto-parallel in an integrable sesquiWalker $\alpha\beta$-geometry, and parallel in the two-sided Walker case. Note that the treatment of $\langle \partial_u \rangle_{\bf R}$ given in Law (2009) (6.1.33) and (6.1.44) is, in regard to the assumptions concerning coordinates and spin frames employed there, consistent with the assumptions for (2.10b) (whereas the coordinates and spin frames employed in Law 2009 (6.1.33) for the treatment of $\langle \partial_v \rangle_{\bf R}$ are not the same as in form (2.10a)), so we will restrict attention here to that form, though the geometric results will be valid generally. The computation of $\alpha_A\alpha^B\nabla_{BB'}\alpha^A = 0$ in the proof of proposition 2.5 above could have been stated as the analogue of Law (2009) (6.2.4), i.e., integrability of $W_{[\lambda]}$ is equivalent (for form (2.10b)) to $\kappa = \sigma = 0$. For oriented Walker coordinates, Law (2009) (5.6) gives $\kappa = 0$ and $\sigma = -b_u/2$, confirming proposition 2.5 in this case. Conditions on the spin coefficients for $\cal D$ to be parallel are stated in Law (2009) (6.1.12) which, together with Law (2009) (5.6), confirm proposition 2.5 for two-sided Walker geometry. The particularly simple forms one obtains in Law (2009) (6.1.33a--c) with $b_u = c_u = 0$ (for form (2.10b)) are consistent with the fact that $\cal D$ is parallel in two-sided Walker geometry (in particular, connecting vector fields between the null geodesics of $\cal D$  within $\cal H$, when expressed in terms of the null tetrad, have constant components of $m^a$ and $\tilde m^a$).

The assumptions underlying form (2.10b) are also consistent with the assumptions employed in the description of null distributions of type III in Law (2009) \S 3. In particular, one confirms that Law (2009) (6.3.9b) is consistent with the curvature results for integrable sesquiWalker $\alpha\beta$-geometries and (6.3.9c) with the curvature results for two-sided Walker geometries. Note that the assertion there that $\Psi_2 = S/12 = \tilde\Psi_2$, where $S$ is the Ricci scalar curvature, is confirmed, for form (2.10b), by Law \& Matsushita (2008) (2.20) and (2.25) ($S = a_{uu} + b_{vv} + 2c_{uv} = a_{uu} + b_{vv}$, when $c_u = 0$). Note that the results in the case of form (2.10b) are entirely consistent with Law (2009) (6.3.13).

Moreover, one also observes that the result of 2.5 above for two-sided Walker geometry is a special case of Walker's (1950b) canonical form for a (four-dimensional) neutral metric admitting a parallel one-dimensional null distribution (equivalently, a parallel three-dimensional null distribution, i.e., of type III in the sense of Law 2009), see Law (2009) (6.1.12) and (6.3.4).

In an integrable sesquiWalker $\alpha\beta$-geometry which is {\sl not} two-sided Walker, Lemma 2.3 entails that ${\cal D} = \langle {\tilde S}^a \rangle_{\bf R}$ (the latter is of course not defined in two-sided Walker geometry). The analogue of $\langle \tilde S^a \rangle_{\bf R}$, i.e., $\langle S^a \rangle_{\bf R}$ in an $\alpha$-geometry which is not Walker, was studied in Law (2009). This coincidence does not appear to be illuminating, however, as the results obtained in  Law (2009) would mainly concern integrability conditions for $\langle \tilde S^a \rangle_{\bf R}^\perp = {\cal H}$, and we have already assumed integrability for $\cal H$ to obtain the coincidence.

We now turn to curvature restrictions in the context of $\alpha\beta$-geometries. Law (2009) (6.2.45) showed that Ricci-null Walker geometries $(M,g,[\pi^{A'}])$, i.e., a Walker geometry for which $[\pi^{A'}]$ is a multiple RPS, take a special form, viz., for Walker coordinates $(u,v,x,y)$ there are functions $\vartheta(u,v,x,y)$, $F(u,x,y)$ and $G(v,x,y)$ satisfying only $F_{uu} = G_{vv} =:h(x,y)$ in terms of which $W$ in (2.8) takes the form
$$W = -2\pmatrix{\vartheta_{vv}&-\vartheta_{uv}\cr -\vartheta_{uv}&\vartheta_{uu}\cr} + \pmatrix{F&0\cr 0&G\cr}.$$
Law (2009) (6.2.47--57) expresses the curvature in terms of $\vartheta$, $F$, and $G$. In particular, the Ricci scalar curvature $S = 2h(x,y)$.
\vskip 24pt
\noindent {\bf 2.6 Proposition}\hfil\break
Let $(M,g,[\pi^{A'}],[\lambda^A])$ be a sesquiWalker $\alpha\beta$-geometry for which $[\pi^{A'}]$ is a multiple RPS. Then, $[\lambda^A]$ is a multiple WPS. Moreover, for any LSR $\lambda^A$ of $[\lambda^A]$, $\phi_{A'B'} := \Phi_{A'B'AB}\lambda^A\lambda^B$ satisfies
$$\lambda^A\nabla^{A'}_A\phi_{A'B'} = 2\zeta^{A'}\phi_{A'B'},$$
where $\zeta_{A'}$ is defined in (2.5b), whence there is an LSR $\chi^A$ of $[\lambda^A]$ (with the freedom to scale by functions constant on $\beta$-surfaces) such that $\chi^A\nabla^{A'}_A\phi_{A'B'} = 0$.

Proof. For specificity, choose oriented Walker coordinates $(u,v,x,y)$ satisfying (2.10a). By 2.5, $a_v = 0$. From Law (2009) (6.2.45d), $a_v = 0\ \Leftrightarrow\ \vartheta_{vvv} = 0$, whence, by Law (2009) (6.2.53), $\Psi_3 = \Psi_4 = 0$, i.e., $[\lambda^A]$ is a multiple WPS. Since $[\lambda^A]$ also satisfies (2.3), by the Generalized Goldberg-Sachs Theorem (GGST) (see e.g., Law 2009, (6.2.17))
$$0 = \lambda^A\lambda^B\lambda^C\nabla^D_{B'}\Psi_{ABCD} = \lambda^A\lambda^B\lambda^C\nabla^{A'}_A\Phi_{BCA'B'},$$
where the second equality follows by a spinor Bianchi identity. Hence
$$\lambda^A\nabla^{A'}_A\phi_{A'B'} = \lambda^A\lambda^B\lambda^C\nabla^{A'}_A\Phi_{BCA'B'} + (\lambda^B\lambda^A\nabla^{A'}_A\lambda^C + \lambda^C\lambda^A\nabla^{A'}_A\lambda^B)\Phi_{BCA'B'} = 2\phi_{A'B'}\zeta^{A'}.$$
Since $[\lambda^A]$ is a solution of (2.3) and a multiple WPS, by Law (2009) (6.2.32), there is an LSR $\chi^A$ of $[\lambda^A]$ for which $\zeta^{A'} = 0$.\bull
\vskip 24pt
The condition $\vartheta_{vvv} = 0$, i.e., $\vartheta$ quadratic in $v$ with coefficients functions of $u$, $x$ and $y$, can be exploited to generate Ricci-null Walker geometries on the chart $(u,v,x,y)$ for which the Walker spin frame element $\beta^A$ defines via $[\beta^A]$ an integrable $\beta$-distribution, i.e., a sesquiWalker $\alpha\beta$-geometry satisfying Proposition 2.6.

One obtains a further specialization of 2.6 if $(M,g,[\pi^{A'}],[\lambda^A])$ is two-sided Walker. The condition $c_v = 0$ is equivalent, by Law (2009) (6.2.45d), to $\vartheta_{uvv} = 0$, i.e., the leading coefficient in the quadratic expression of $\vartheta$ as a function of $v$ is independent of $u$. One then notes from Law (2009) (6.2.53) that  $\Psi_2 = h/6 = S/12 = \tilde\Psi_2$ and from Law (2009) (6.2.57) that $A_{11} = 0$, i.e., $[\lambda^A]$ is a RPS, these results confirming those obtained above for any two-sided Walker geometry.

In a self-dual (SD) Walker geometry $(M,g,[\pi^{A'}])$, the absence of ASD Weyl curvature means that at each point $p$ and for each $\beta$-plane $W \subseteq T_pM$, there exists a $\beta$-surface $\cal W$ such that $T_p{\cal W} = W$ (see, for example, LeBrun \& Mason 2007, \S 3). In fact, one can construct, locally, integrable $\beta$-distributions. Let $[\lambda^A]$ be a (local) spinor field defining an integrable $\beta$-distribution so that $(M,g,[\pi^{A'}],[\lambda^A])$ is SD sesquiWalker. D\'{\i}az-Ramos et al. (2006) and Davidov \& Mu\v skarov (2006) provided  a local characterization of SD Walker metrics. With respect to Walker coordinates $(u,v,x,y)$, a Walker metric is SD iff the metric components $a$, $b$, and $c$ in (2.8) take the form:
$$\vcenter{\openup1\jot \halign{$\hfil#$&&${}#\hfil$&\qquad$\hfil#$\cr
a &= Au^3 + Bu^2v + Cu^2 + 2Duv + Eu + Fv + G;\cr
b &= Bv^3 + Auv^2  + Kv^2 + 2Luv + Mu + Nv + H;\cr
c &= Au^2v +  Buv^2 + Lu^2 + Dv^2 + {C+K \over 2}uv + Pu + Qv + T;\cr}}\eqno(2.15)$$
where the coefficients are arbitrary functions of $x$ and $y$. Assuming the geometry is in fact integrable sesquiWalker, choosing Walker oriented coordinates satisfying (2.10a), then $a_v = 0$, i.e., $B=D=F=0$.  If one further supposes the geometry is two-sided Walker, then $c_v = 0$ too, i.e., altogether one has $A = B = D = F = Q = 0$ and $C = -K$, resulting in the simpler equations:
$$\vcenter{\openup1\jot \halign{$\hfil#$&&${}#\hfil$&\qquad$\hfil#$\cr
a &= Cu^2 + Eu + + G;\cr
b &= -Cv^2 + 2Luv + Mu + Nv + H;\cr
c &= Lu^2 + Pu + T.\cr}}\eqno(2.16)$$
Now suppose $(M,g,[\pi^{A'}])$ is a SD Ricci-null Walker geometry. Once again, due to the self duality, one can introduce an integrable $\beta$-distribution and consider Ricci-null, SD integrable sesquiWalker $\alpha\beta$-geometries $(M,g,[\pi^{A'}],[\lambda^A])$ using the description in Law (2009), (6.2.45). The SD condition, by Law (2009) (6.2.53) is equivalent to
$$\vartheta_{vvvv} = \vartheta_{uvvv} = \vartheta_{uuuv} = \vartheta_{uuuu} = 0 \hbox{  and  } \vartheta_{uuvv} = {h \over 6}.\eqno(2.17)$$
Supposing the oriented Walker coordinates $(u,v,x,y)$ have been chosen to satisfy (2.10a), then $a_v = 0$, i.e., $\vartheta_{vvv} = 0$, which is compatible with (2.17). A simplification is obtained by assuming the Ricci scalar curvature vanishes ($S = 0$), which is equivalent to $h=0$. Then, all fourth order partial derivatives in $u$ and $v$ of $\vartheta$ vanish and $\vartheta$ is a cubic polynomial in $u$ and $v$, with coefficients functions of $x$ and $y$; the condition $\vartheta_{vvv} = 0$ excludes a term in $v^3$.

If one considers a SD, two-sided Walker geometry $(M,g,[\pi^{A'}],[\lambda^A])$ for which $[\pi^{A'}]$ is a multiple RPS, then as noted previously for two-sided Walker geometries, $\Psi_2 = S/12 = \tilde\Psi_2$. Hence, self duality entails $0 = \Psi_2 = S/12 = \tilde\Psi_2$, and $[\pi^{A'}]$ is a WPS of multiplicity at least three. If the oriented Walker coordinates $(u,v,x,y)$ are chosen to satisfy (2.10a), then $a_v = c_v = 0$, i.e., $\vartheta_{vvv} = \vartheta_{uvv} = 0$. Together with the vanishing of all the fourth order partial derivatives in $u$ and $v$, $\vartheta$ must therefore take the form
$$\vartheta = K_2u^2v + K_4u^3 + K_5u^2 + K_6uv + K_7v^2 + K_8u + K_9v + K_{10},\eqno(2.18)$$
where the coefficients are functions of $x$ and $y$. One special case is a left-flat Walker geometry $(M,g,[\pi^{A'}])$, i.e., a Walker geometry for which the SD Weyl curvature is the only nontrivial curvature. Such curvature permits the construction, locally, of parallel unprimed spin frames, whence, locally, of parallel $\beta$-distributions, so left-flat Walker geometries are, locally, automatically two-sided Walker, and Ricci null with respect to $[\pi^{A'}]$. As $S = 0$, from Law (2009) (6.2.58), $F(u,x,y) = uf(x,y)$ and $G(v,x,y) = vg(x,y)$, for some functions $f$ and $g$. Hence, locally, such geometries are characterized by the expression (2.18) satisfying Law (2009) (6.2.63). Substitution of (2.18) into Law (2009) (6.2.63) yields two constraints:
$$(8K_2 - g)_x = f_y \hskip 1in (3K_4)_x + (K_2)_y - 4(K_2)^2 + {g \over 2}K_2 - {f \over 2}(3K_4) = 0.$$
The first implies the existence of a function $X(x,y)$ satisfying $X_x = f$ and $X_y = 8K_2 - g$. Upon substituting these expressions into the second constraint, one obtains
$$(3K_4)_x + (K_2)_y = {X_x \over 2}(3K_4) + {X_y \over 2}K_2,$$
which is of the form of Chudecki \& Przanowski (2008b) (5.35) and therefore has solution $3K_4 = Y_y\exp(X/2)$ and $K_2 = -Y_x\exp(X/2)$, for some function $Y(x,y)$. Substituting into Law (2009) (6.2.59) yields
$$a = uX_x - 4K_7 \qquad c = -4uY_x\exp(X/2)+2K_6 \qquad b = -4(uY_y + vY_x)\exp(X/2) - vX_y - 4K_5,\eqno(2.19)$$
which is essentially Chudecki \& Przanowski (2008b) (5.37).
\vskip 24pt
\noindent{\section 3. Algebraically Special AlphaBeta-Geometries}
\vskip 12pt
Let $(M,h,[\pi^{A'}],[\lambda^A])$ be an $\alpha\beta$-geometry. As noted in \S 2, $[\pi^{A'}]$ and $[\lambda^A]$ are each WPSs. In this section we assume that $[\pi^{A'}]$ is a multiple WPS and say that $(M,h,[\pi^{A'}])$ is a real algebraically special (AS)$\alpha$-geometry. By Law \& Matsushita (2011), 3.15, such is locally conformal to a Walker geometry, i.e., each point $p \in M$ has a neighbourhood $U$ such that on $U$, $h = \Omega^2g$, for some smooth function $\Omega:U \to {\bf R}^+$ and some metric $g$ on $U$ such that $(U,g,[\pi^{A'}])$ is Walker. Since $\alpha$- and $\beta$-distributions are conformally invariant structures and integrability of distributions is a differential-topological condition, $(U,g,[\pi^{A'}],[\lambda^A])$ is a sesquiWalker $\alpha\beta$-geometry. Moreover, the integrability of $\cal H$ is also unaffected by the conformal rescaling, whence one can evaluate that condition on $U$ with respect to $(U,g,[\pi^{A'}],[\lambda^A])$ using lemma 2.3. Granted integrability of $\cal H$, one can then deduce properties of the local geometry of $(M,h,[\pi^{A'}],[\lambda^A])$ using conformal rescaling formulae as in Law \& Matsushita (2011) and the known local geometric properties of integrable sesquiWalker $\alpha\beta$-geometry.

Now suppose that $[\lambda^A]$ is a {\sl multiple} WPS too. As this property is conformally invariant, it remains valid for $(U,g,[\pi^{A'}],[\lambda^A])$, whence $(U,g,[\lambda^A])$ is a real AS$\beta$-geometry and it is locally conformal to a geometry which is Walker for the $\beta$-distribution, i.e., each point $p \in U$ has a neighbourhood $V$ on which $g = \chi^2k$, for some smooth $\chi:V \to {\bf R}^+$, where $(V,k,[\lambda^A])$ is Walker. We will refer to $(V,k,[\pi^{A'}],[\lambda^A])$ as a sesquiWalker $\beta\alpha$-geometry. 

A natural question to ask is when $(V,k,[\pi^{A'}],[\lambda^A])$ can be constructed so as to be two-sided Walker, i.e., in effect, when is $(M,h,[\pi^{A'}],[\lambda^A])$ locally conformal to two-sided Walker geometry? Note that a necessary condition is that $\cal H$ must be integrable as it is so for two-sided Walker geometry and is a differential-topological condition.

If one begins with an $\alpha\beta$-geometry $(M,h,[\pi^{A'}],[\lambda^A])$, and conformally rescales by $\Omega^2$, from (2.2) and (2.6) one requires
$$\pi_{A'}\hat\nabla_b\pi^{A'} = 0\qquad\hbox{and}\qquad \lambda_A\hat\nabla_b\lambda^A = 0,$$
for the rescaled geometry to be two-sided Walker, i.e.,
$$0 = S_b + \pi_{B'}\Upsilon_{BX'}\pi^{X'} \hskip 1in 0 = \tilde S_b + \lambda_B\Upsilon_{XB'}\lambda^X,$$
where $\Upsilon_a = \nabla_a\ln\Omega$ (see, e.g., Law \& Matsushita 2011), i.e., one must solve
$$\pi^{B'}\nabla_{BB'}f = \omega_B \hskip 1in \lambda^B\nabla_{BB'}f = \kappa_{B'},\eqno(3.1)$$
for some $f$. Taking components gives four equations but two are equivalent under the assumption that $\cal H$ is integrable.

In fact, from Law \& Matsushita (2011), a necessary and sufficient condition to solve the first equation of (3.1) is that $[\pi^{A'}]$ is a multiple WPS, and to solve the second equation is that $[\lambda^A]$ is a WPS. Moreover, one obviously has the necessary condition that
$$\pi^{A'}\kappa_{A'} = \lambda^A\pi^{A'}\nabla_{AA'}f = \lambda^A\omega_A,$$
which is (2.7), i.e., integrability of $\cal H$. A natural question is whether these three necessary conditions, integrability of $\cal H$ and each of $[\pi^{A'}]$ and $[\lambda^A]$ {\sl multiple} WPSs, are also sufficient to simultaneously solve the pair of equations in (3.1)? 

To investigate this question, note that if $(U,\Omega^{-2}h,[\pi^{A'}],[\lambda^A])$ is two-sided Walker for some smooth function $\Omega:U \to {\bf R}^+$, then, by Law \& Matsushita (2011), 3.1, for any smooth function $\chi:U \to {\bf R}^+$ which is constant on $\alpha$-surfaces, the $\alpha$-distribution is still parallel in $(U,\chi^2\Omega^{-2}h,[\pi^{A'}])$. Hence, if it is possible to locally conformally rescale $h$ to be two-sided Walker, then it is possible to first rescale $h$ to be Walker for $[\pi^{A'}]$ and then to locally rescale that Walker metric to be two-sided Walker. We take this approach to study our question so as to utilize the Walker coordinates introduced in \S 2. 

Since $(M,h,[\pi^{A'}],[\lambda^A])$ is an AS$\alpha$-geometry, as noted in the opening paragraph of this section, each point $p \in M$ has a neighbourhood $U$ such that $(U,\Omega^{-2}h,[\pi^{A'}],[\lambda^A])$ is a sesquiWalker $\alpha\beta$-geometry for some $\Omega:U \to {\bf R}^+$. Can we choose, on a possibly smaller neighbourhood $V$ of $p$, a smooth $\chi:V \to {\bf R}^+$ such that $(V,\chi^2\Omega^{-2}h,[\pi^{A'}],[\lambda^A])$ is two-sided Walker? One would require that $\chi$ solve $\pi^{B'}\nabla_{BB'}\chi = 0$ and $\lambda^B\nabla_{BB'}\left[\ln(\chi^{-1})\right] = \kappa_{B'}$, where these equations are the appropriate analogues of Law \& Matsushita (2011) (3.15.1), i.e., the equations of (3.1) above in the present circumstances. For simplicity, we may write these equations as
$$\pi^{B'}\nabla_{BB'}f = 0, \hskip .75in \lambda^B\nabla_{BB'}f = \kappa_{B'}, \hskip .75in f := \ln\left(\chi^{-1}\right).\eqno(3.2)$$
Since we assume $\cal H$ is integrable, $(U,\Omega^{-2}h,[\pi^{A'}],[\lambda^A])$ is integrable sesquiWalker and one can choose oriented Walker coordinates $(u,v,x,y)$ satisfying (2.10a). Let $\pi^{A'}$ and $\lambda^A$ be the LSRs of $[\pi^{A'}]$ and $[\lambda^A]$, respectively, determined by (2.10a) (determined up to a common sign). In the Walker geometry $(U,\Omega^{-2}h,[\pi^{A'}])$, $\omega_A = 0$ and integrability of $\cal H$ is the condition $\pi^{A'}\kappa_{A'} =0$. Indeed, from the proof of Proposition 2.5, in $(U,\Omega^{-2}h,[\pi^{A'}])$: $\lambda_A\kappa_{A'} = \lambda_B\nabla_{AA'}\lambda^B = \beta_B\nabla_{AA'}\beta^B = -(c_v/2)\tilde m_a$ ($a_v = 0$ since the $\beta$-distribution is integrable), i.e.,
$$\kappa_{B'} = -{c_v \over 2}\pi_{B'}.\eqno(3.3)$$
Thus, as must be the case, $\kappa_{B'} = 0$ iff $c_v = 0$, i.e., iff $(U,\Omega^{-2}h,[\pi^{A'}],[\lambda^A])$ is already two-sided Walker.

The equations (3.2) are a system of PDEs on the integral surfaces of $\cal H$. The coordinates $(v,u,x,y)$ are Frobenius coordinates for the nested distributions ${\cal D} \leq Z_{[\pi]} \leq {\cal H}$. Taking components, by (2.10a--11a), (3.2) are equivalent to
$$X_1f := \partial_vf = 0, \hskip .5in X_2f := \partial_uf = 0, \hskip .5in X_3f := -\lambda^B\xi^{B'}\nabla_{BB'}f = -\xi^{B'}\kappa_{B'} = {c_v \over 2}.$$
By Law \& Matsushita (2008), (2.11),
$$X_3f = -n^b\nabla_bf = \left({a \over 2}\partial_u + {c \over 2}\partial_v - \partial_x\right)f = {c_v \over 2}.$$
Writing $[X_i,X_j] =: \Phi^k_{ij}X_k$, one computes: $[X_1,X_2] = 0$, $[X_1,X_3] = (a_v/2)X_2 + (c_v/2)X_1 = (c_v/2)X_1$, $[X_2,X_3] = (a_u/2)X_2 + (c_u/2)X_1$. Hence, the only nonzero $\Phi^k_{ij}$ are
$$\Phi^1_{13} = {c_v \over 2}, \hskip .75in \Phi^1_{23} = {c_u \over 2}, \hskip .75in \Phi^2_{23} = {a_u \over 2}.$$ 
With $\phi_1 = \phi_2 = 0$ and $\phi_3 = c_v/2$, the integrability conditions for $X_if = \phi_i$, $i=$1--3, are $X_i\phi_j - X_j\phi_i = \Phi^k_{ij}\phi_k$. The non vacuous conditions are: $X_1\phi_3 - X_3\phi_1 = \Phi^k_{13}\phi_k$, i.e., $(c_{vv}/2) = 0$; $X_2\phi_3 - X_3\phi_2 = \Phi^k_{23}\phi_k$, i.e., $(c_{vu}/2) = 0$. Thus, the integrability conditions for (3.2), in these Walker coordinates, are
$$c_{uv} = c_{vv} = 0.$$
As $a_v = 0$ already, from Law \& Matsushita (2008), (2.25), $[\lambda^A] = [\beta^A]$ is a multiple WPS iff $c_{vv} = 0$. Thus, granted our assumptions, the remaining obstruction to solving (3.2) in terms of the Walker coordinates satisfying (2.10a) is the single condition
$$c_{uv} = 0,\eqno(3.4)$$
clearly weaker than $c_v = 0$, the condition here for $(U,\Omega^{-2}h,[\pi^{A'}])$ to be two-sided Walker. Exactly the same result is obtained utilizing the coordinate form (2.10b). The condition (3.4) indicates that the multiplicity of each of $[\pi^{A'}]$ and $[\lambda^A]$ as WPSs together with the integrability of $\cal H$ are necessary but not sufficient conditions for $(M,h,[\pi^{A'}],[\lambda^A])$ to be locally conformally two-sided Walker. A geometric characterization of this obstruction is desirable. From Law \& Matsushita (2008) (2.25), (3.4) is equivalent to
$$\Psi_2 = {S \over 12},\eqno(3.5)$$
for the Walker geometry $(U,\Omega^{-2}h,[\pi^{A'}])$. Given that $[\beta^A]$ is here a multiple WPS ($\Psi_4 = \Psi_3 = 0$) and that the freedom in the spin frame $\{\alpha^A,\beta^A\}$ is $\beta^A \mapsto \lambda\beta^A$, $\alpha^A \mapsto \lambda^{-1}\alpha^A + \mu\beta^A$, then $\Psi_2 = \Psi_{ABCD}\alpha^A\alpha^B\beta^C\beta^D$ is in fact a geometric quantity here. Thus, under the stated constraints, (3.5) is a geometric condition, albeit expressed in terms of the Walker geometry $(U,\Omega^{-2}h,[\pi^{A'}])$ rather than the $\alpha\beta$-geometry $(M,h,[\pi^{A'}],[\lambda^A])$. (Note that while (3.4) is equivalent to (3.5) for arbitrary Walker coordinates, $\Psi_2$ is only geometric in nature in the circumstances stated above.)

Under any conformal rescaling $g \mapsto \chi^2g$, exploiting the spin frames of Law \& Matsushita (2011) (3.7) for the rescaled metric, one has
$$\hat \Psi_2 = \chi^{-2}\Psi_2 \hskip 1in \hat S = \chi^{-2}[S - 6\chi^{-1}\Square\chi]\eqno(3.6)$$
(see, for example, Law \& Matsushita 2011, (2.7)). In Walker geometry, by Law \& Matsushita (2008), (3.9),
$$\Square\chi = -a\chi_{uu} - 2c\chi_{uv} - b\chi_{vv} + 2\chi_{ux} + 2\chi_{vy} - (a_u+c_v)\chi_u - (b_v+c_u)\chi_v.\eqno(3.7)$$
Note that if $\chi$ is constant on $\alpha$-surfaces, i.e., is a function of $(x,y)$ only, then $\Square\chi = 0$ and under such conformal rescalings (3.5) is invariant. Thus, choosing a $\chi$ to solve the first equation of (3.2), the obstruction to solving the second equation of (3.2) is again (3.5) within the rescaled geometry. As noted previously, (3.5) follows when $\cal H$ is parallel and thus is indeed a necessary condition if $(V,\chi^2\Omega^{-2}h,[\pi^{A'}],[\lambda^A])$ is to be two-sided Walker for some choice of $\chi$. Thus, we have proved the following result.
\vskip 24pt
\noindent {\bf 3.1 Proposition}\hfil\break
An $\alpha\beta$-geometry $(M,h,[\pi^{A'}],[\lambda^A])$ is locally conformally two-sided Walker iff $\cal H$ is integrable, each of $[\pi^{A'}]$ and $[\lambda^A]$ are multiple WPSs, and (3.5) holds in every locally conformally-related Walker geometry $(U,\Omega^{-2}h,[\pi^{A'}])$.
\vskip 24pt
We end with some observations on an explicit example presented by Chudecki \& Przanowski (2008a), (4.23), which we write in the form:
$$\eqalignno{(h_{\bf ab}) &= v^{-2}(g_{\bf ab})&(3.8)\cr
&= v^{-2}\left[2(dudx + dvdy) + \left({e^{4F}u^4 \over 3v^2} + 4uF_x\right)dx^2 + 2\left({2e^{4F}u^3 \over 3v} + 2uF_y\right)dxdy + \left(e^{4F}u^2 + 2vF_y\right)dy^2\right],\cr}$$
where $F$ is a function of $(x,y)$ only. It is clear that the metric $(g_{\bf ab})$ is of the form (2.8), i.e., is Walker, with $(u,v,x,y)$ Walker coordinates. The metric $(h_{\bf ab})$ is obviously conformally Walker. Chudecki \& Przanowski (2008a) constructed this metric as a solution of the hyperheavenly equation and thus $h$ is Einstein; moreover, they constructed it so that each Weyl spinor is of type $\{4\}$ (i.e., null). We shall denote the multiple WPSs by $[\pi^{A'}]$ and $[\lambda^A]$; since they are the only WPSs, $[\pi^{A'}]$ must be the projective spinor field defining the integrable $\alpha$-distribution. It follows that $h$ is not itself Walker for this $\alpha$-distribution $Z_{[\pi]}$ since the conformal factor $\chi := v^{-1}$ is not constant on $\alpha$-surfaces; there are no other integrable $\alpha$-distributions as there is only one WPS for $\tilde\Psi_{A'B'C'D'}$. There is also at most one integrable $\beta$-distribution. Indeed, since the metric $h$ is Einstein, a spinor Bianchi identity gives $\nabla^A_{B'}\Psi_{ABCD} = 0$, whence the GGST (Law 2009 (6.2.19)) ensures that $[\lambda^A]$ is in fact a solution of (2.3) and does define an integrable $\beta$-distribution $W_{[\lambda]}$. 

Since $[\pi^{A'}]$ is of multiplicity four, the Ricci scalar curvature $S$ of the Walker metric $g$ must vanish (Law \& Matsushita (2008), 2.5 or 2.6). The Ricci scalar curvature $\hat S$ of the metric $h$ is then
$$\hat S = \chi^{-2}[S - 6\chi^{-1}\Square\chi] = -6\chi^{-3}\Square\chi.$$
From (3.7) one computes that $\Square\chi = 0$, whence $\hat S = 0$. From Law \& Matsushita (2008) A1.8 and (2.33), however, one easily checks that $[\pi^{A'}]$ is not a multiple RPS for the Walker metric $g$ (the $\theta$, $\mu$, and $\nu$ in (2.33) are each nonzero); in particular, the Walker metric $g$ is not Einstein. 

Using the given Walker coordinates $(u,v,x,y)$ for $g$, and their associated Walker spin frames, one can readily calculate the ASD Weyl curvature spinor from Law \& Matsushita (2008) (2.25--26), obtaining:
$$\Psi_{ABCD} = e^{4F}\left({u^4 \over v^4}\alpha_A\alpha_B\alpha_C\alpha_D + 4{u^3 \over v^3}\alpha_{(A}\alpha_B\alpha_C\beta_{D)} + 6{u^2 \over v^2}\alpha_{(A}\alpha_B\beta_C\beta_{D)} + 4{u \over v}\alpha_{(A}\beta_B\beta_C\beta_{D)} + \beta_A\beta_B\beta_C\beta_D\right).$$
To determine the WPS of $\Psi_{ABCD}$, one need only equate this expression with $\lambda_A\lambda_B\lambda_C\lambda_D$. One finds $[\lambda^A] = [u\alpha^A + v\beta^A]$. Using (2.12) and the covariant derivatives of the Walker spin frames (Law 2009, (5.8)) one readily computes:
$$\lambda_A\nabla_b\lambda^A = \left({e^{4F}u^4 \over 6v}\right)\ell_b + \left({e^{4F}u^3 \over 6}\right)\tilde m_b - vn_b - um_b = \lambda_B\left[{e^{4F}u^3 \over 6v}\pi_{B'} - \xi_{B'}\right] \not= 0 \hskip .5in \lambda_A\lambda^B\nabla_{BB'}\lambda^A = 0.$$
Hence the Walker metric $g$ is indeed sesquiWalker, but not two-sided Walker. Moreover, from (2.5b), $\kappa_{B'} = (e^{4F}u^3/6v)\pi_{B'} - \xi_{B'}$, so $\kappa^{B'}\pi_{B'} = -1$ whence, by lemma 2.3, $\cal H$ is not integrable. Thus, the metric $h$, though an $\alpha\beta$-metric which is algebraically special for both WPSs, is not locally conformally two-sided Walker. Note, however, that the condition (3.5) is valid for the metric $g$ (both sides of the equation vanish), so the non-integrability of $\cal H$ is the only obstruction to $h$ being locally conformally two-sided Walker.
\vskip 24pt
\noindent {\section References}
\vskip 12pt
\frenchspacing
\hangindent=20pt \hangafter=1
\noindent Brozos-V\'azquez, M., Garc\'{\i}a-R\'{\i}o, E., Gilkey, P., Nik\v cevi\'c, S. \& V\'azquez-Lorenzo, R. 2009 {\sl The Geometry of Walker Manifolds}. Morgan \& Claypool. 
\vskip 1pt
\hangindent=20pt \hangafter=1
\noindent Boyer, C. P., Finley III, J. D. \& Pleba\'nski, J. F. 1980 Complex General Relativity, $\cal H$ and ${\cal HH}$ Spaces-A Survey of One Approach, in {\sl General Relativity and Gravitation: One Hundred Years After the Birth of Albert Einstein}, Vol. 2, A. Held (ed.), Plenum Press, New York, NY, 241--281.
\vskip 1pt
\hangindent=20pt \hangafter=1
\noindent Chudecki, A. and Przanowski, M. 2008a A simple example of type-$[N] \otimes [N] {\cal HH}$-spaces admitting twisting null geodesic congruence. {\sl Classical Quantum Gravity}, {\bf 25}, 055010 (13pp).
\vskip 1pt
\hangindent=20pt \hangafter=1
\noindent Chudecki, A. and Przanowski, M. 2008b From hyperheavenly spaces to Walker spaces and Osserman spaces. I. {\sl Classical Quantum Gravity}, {\bf 25} (145010) (18 pp).
\vskip 1pt
\hangindent=20pt \hangafter=1
\noindent Davidov, J. \& Mu\v skarov, O. 2006 Self-dual Walker metrics with two-step nilpotent Ricci operator. {\sl J. Geom. Phys.} {\bf 57}, 157--165.
\vskip 1pt
\hangindent=20pt \hangafter=1
\noindent D\'{\i}az-Ramos, J. C., Garc\'\i a-R\'\i o, E. \& V\'azquez-Lorenzo, R. 2006 New Examples of Osserman metrics with nondiagonalizable Jacobi operators. {\sl Differential Geometry and Its Applications\/} {\bf 24}, 433--442.
\vskip 1pt
\hangindent=20pt \hangafter=1
\noindent Finley III, J. D. \& Pleba\'nski, J. F. 1976 The intrinsic spinorial structure of hyperheavens. {\sl Journal of Mathematical Physics} {\bf 17}, 2207--2214.
\vskip 1pt
\hangindent=20pt \hangafter=1
\noindent Gover, A. R., Hill, C. D. and Nurowski, P. 2011 Sharp version of the Goldberg-Sachs theorem. {\sl Annali di Matematica} {\bf 190}, 295--340.
\vskip 1pt
\hangindent=20pt \hangafter=1
\noindent Hansen, R. O., Newman, E. T., Penrose, R. and Tod, K. P. 1978 The metric and curvature properties of ${\cal H}$-space. {\sl Proceedings of the Royal Society London A} {\bf 363}, 445--468.
\vskip 1pt
\hangindent=20pt \hangafter=1
\noindent Ko, M., Ludvigsen, M., Newman, E. T. \& Tod, K. P. 1981 The theory of ${\cal H}$-space. {\sl Physics Reports} {\bf 71}, 51--139.
\vskip 1pt
\hangindent=20pt \hangafter=1
\noindent Law, P. R. 2006 Classification of the Weyl curvature spinors of neutral metrics in four dimensions. {\sl J. Geo. Phys.} {\bf 56}, 2093--2108.
\vskip 1pt
\hangindent=20pt \hangafter=1
\noindent Law, P. R. 2009 Spin Coefficients for Four-Dimensional Neutral Metrics, and Null Geometry. {\sl Journal of Geometry and Physics} {\bf 59 (8)}, 1087--1126; arXiv:0802.1761v2 [math.DG] (26 Aug 2009).
\vskip 1pt
\hangindent=20pt \hangafter=1
\noindent Law, P. R. \& Matsushita, Y. 2008 A Spinor Approach to Walker Geometry. {\sl Communications in Mathematical Physics}, {\bf 282}, 577-623; arXiv:math/0612804v4 [math.DG] 7 Apr 2009.
\vskip 1pt
\hangindent=20pt \hangafter=1
\noindent Law, P. R. \& Matsushita, Y. 2011 Algebraically Special, Real Alpha-Geometries. {\sl Journal of Geometry and Physics} {\bf 61}, 2064--2080. arXiv:0808.2082v3 [math.DG] (9 Sep 2011).
\vskip 1pt
\hangindent=20pt \hangafter=1
\noindent LeBrun, C. and Mason, L. J. 2007 Nonlinear Gravitons, Null Geodesics, and Holomorphic Disks. {\sl Duke Mathematics Journal} {\bf 136}, 205--273.
\vskip 1pt
\hangindent=20pt \hangafter=1
\noindent Matsushita, Y., Haze, S. and Law, P. R. 2007 Almost K\"ahler-Einstein Structures on 8-Dimensional Walker Manifolds. {\sl Monatsh. Math.} {\bf 150}, 41--48.
\vskip 1pt
\hangindent=20pt \hangafter=1
\noindent Newman, E. T. 1976 Heaven and its properties. {\sl General Relativity and Gravitation} {\bf 7}, 107-111.
\vskip 1pt
\hangindent=20pt \hangafter=1
\noindent Newman, E. T. and Tod, K. P. 1981 A note on left-flat space-times. {\sl Journal of Mathematical Physics} {\bf 21}, 874--877.
\vskip 1pt
\hangindent=20pt \hangafter=1
\noindent Penrose, R. Nonlinear Gravitons and Curved Twistor Theory {\sl General Relativity and Gravitation} {\bf 7}, 31--52.
\vskip 1pt
\hangindent=20pt \hangafter=1
\noindent Penrose, R. 1999 Twistor theory and the Einstein vacuum. {\sl Classical and Quantum Gravity} {\bf 16}, A113--A130.
\vskip 1pt
\hangindent=20pt \hangafter=1
\noindent Penrose, R. and Ward, R. S. 1980 Twistors for Flat and Curved Space-Time; In {\sl General relativity and gravitation, Vol. 2}, Plenum, New York and London, pp. 283--328
\vskip 1pt
\hangindent=20pt \hangafter=1
\noindent Pleba\'nski, J. F. 1975 Some solutions of the complex Einstein equations. {\sl Journal of Mathematical Physics} {\bf 16}, 2399--2402.
\vskip 1pt
\hangindent=20pt \hangafter=1
\noindent Pleba\'nski, J. F. and Robinson, I. 1976 Left-degenerate vacuum metrics. {\sl Physical Review Letters} {\bf 37}, 493--495.
\vskip 1pt
\hangindent=20pt \hangafter=1
\noindent Walker, A. G. 1950a Canonical form for a Riemannian space with a parallel field of null planes. {\sl Quart. J. Math. Oxford(2)\/} {\bf 1}, 69--79.
\vskip 1pt
\hangindent=20pt \hangafter=1
\noindent Walker, A. G. 1950b Canonical Forms (II): Parallel Partially Null Planes. {\sl Quart. J. Math. Oxford(2)\/} {\bf 1}, 147--152.
\vskip 1pt
\hangindent=20pt \hangafter=1
\noindent Ward, R. S. 1980 Self-dual space-times with cosmological constant. {\sl Communications in Mathematical Physics} {\bf 78}, 1--17.

\bye